\documentclass{article}
\usepackage[utf8]{inputenc}
\usepackage{tikz}
\usepackage{amssymb,amsmath,caption,subcaption,graphicx,hyperref}
\usepackage{amsthm}

\newcommand{\aw}{\operatorname{aw}}

\newcommand{\dist}{\operatorname{d}}

\newtheorem{thm}{Theorem}[section]

\newtheorem{prop}[thm]{Proposition}
\newtheorem{cor}[thm]{Corollary}
\newtheorem{lem}[thm]{Lemma}

\newtheorem{conj}[thm]{Conjecture}

\newtheorem{obs}[thm]{Observation}

\usepackage{fix-cm}    
\usepackage{float}    

\tikzset{
mystyle/.style={
  circle,
  inner sep=0pt,
  text width=6mm,
  align=center,
  draw=black,
  fill=white
  }
}
\usetikzlibrary{shapes,snakes}
\usetikzlibrary{decorations.pathreplacing}

\makeatletter
\newcommand\HUGE{\@setfontsize\Huge{50}{60}}
\makeatother  
\newcommand{\bpf}{\begin{proof}}
\newcommand{\epf}{\end{proof}}

\newlength\tindent
\renewcommand{\indent}{\hspace*{\tindent}}

\title{Anti-van der Waerden numbers of Graph Products}
\author{Hunter Rehm, Alex Schulte and Nathan Warnberg}
\date{\today}

\begin{document}

\maketitle
\section*{Abstract}

In this paper, anti-van der Waerden numbers on Cartesian products of graphs are investigated and a conjecture made by Schulte, et al (see \cite{SWY}) is answered.  In particular, the anti-van der Waerden number of the Cartesian product of two graphs has an upper bound of four.  This result is then used to determine the anti-van der Waerden number for any Cartesian product of two paths. 

\section{Introduction}

\indent \indent	The anti-van der Waerden number was first defined in \cite{U}. Many results on arithmetic progressions of $[n]$ and the cyclic groups $\mathbb{Z}_n$ were considered in \cite{DMS} and a function $f(n)$ was established in \cite{BSY} such that $\aw([n],3)=f(n)$ for all $n \in \mathbb{N}$.  Results on colorings of the integers with no rainbow $3$-term arithmetic progressions were also studied in \cite{AF} and \cite{AM}.  Colorings and $3$-term arithmetic progressions have been extended to groups (see \cite{finabgroup}) and graphs (see \cite{SWY}).  The authors in \cite{SWY} were inspired to investigate the anti-van der Waerden number of graphs by extending results on the anti-van der Waerden number of $[n]$ and $\mathbb{Z}_n$ to paths and cycles, respectively.  In particular, they noticed that the set of arithmetic progressions on $[n]$ is isomorphic to the set of non-degenerate arithmetic progressions on $P_n$. Similarly, the set of arithmetic progressions on $\mathbb{Z}_n$ is isomorphic to the set of non-degenerate arithmetic progressions on $C_n$. Therefore, considering the anti-van der Waerden number of $[n]$ or $\mathbb{Z}_n$ is equivalent to studying the anti-van der Waerden number of paths or cycles respectively.  The authors of \cite{SWY} made a conjecture about graph products and this conjecture is proven in this paper.  First, some terminology and notation is introduced. \\

A \emph{graph}, $G$, is a collection of vertices, $V(G)$, and edges, $E(G)$, and will be denoted as $G= (V,E)$.  The edge set $E$ is a set of pairs of vertices that indicate the two vertices are connected.  Thus, if there is an edge connecting vertices $u$ and $v$, then $\{u,v\}$ is an edge or $uv$ is an edge for short.  Graph $H$ is a \emph{subgraph} of $G$ if $V(H) \subseteq V(G)$ and $E(H) \subseteq E(G)$.  An \emph{induced subgraph} $H$ of $G$ is one formed by deleting vertices of $G$ and keeping all possible edges.  For the purposes of this paper all graphs are simple (loop free, undirected, no edge weights, no multiple edges) and connected.  The \emph{distance} between vertex $v$ and $u$ in graph $G$ is denoted $\dist_G{v,u}$, $\dist(uv,)$ will be used when the context is clear.  If $G = (V,E)$ and $H = (V',E')$ then the \emph{Cartesian product}, written $G\square H$, has vertex set $\{(x,y) \,|\, x \in V \text{ and } y\in V'\}$ and $(x,y)$ and $(x',y')$ are adjacent in $G\square H$ if either $x = x'$ and $yy'\in E'$ or $y = y'$ and $xx'\in E$.  In this paper, $P_n$ denotes the path graph on $n$ vertices. 

The vertex set of $P_m\square P_n$ is given by $\{v_{i,j} \, | \, 1\le i \le m \text{ and } 1 \le j \le n\}$.  Further, $v_{i,j}$ can be found at the intersection of the $i$th row and $j$th column of $P_m\square P_n$.  This convention allows for the computation of distances in grid graphs based on the subscripts of the vertices.  In particular, if $v_{i,j}$ and $v_{\ell,k}$ are in $P_m\square P_n$ then $\dist(v_{i,j}, v_{\ell,k}) = |i-\ell| + |j-k|$.

\indent A \emph{$k$-term arithmetic progression of a graph} $G$, $k$-AP, is a subset of $k$ vertices of $G$ of the form $\{v_1, v_2, \ldots, v_k \}$, where $\dist(v_i,v_{i+1})=d$ for all $1\leq i<k$.  A $k$-term arithmetic progression is \emph{degenerate} if $v_i=v_j$ for any $i\neq j$.

\indent  An \emph{exact $r$-coloring of a graph} $G$ is a surjective function $c:V(G) \to \{1,2, \dots, r\}$.  A set of vertices $S \subseteq V(G)$ is \emph{rainbow} under coloring $c$, if for any $v_i, v_j \in S$, $c(v_i) \neq c(v_j)$ when $v_i \neq v_j$.  Note that degenerate $k$-APs will not be rainbow.  Given a set of vertices $S\subseteq V(G)$, $c(S) = \{c(v) | v \in S\}$, is the set of colors used on the vertices of $S$.

	The \emph{anti-van der Waerden number of a graph} $G$, denoted by $\aw(G,k)$, is the least positive integer $r$ such that every exact $r$-coloring of $G$ contains a rainbow $k$-AP. If $G$ has $n$ vertices and no coloring of $G$ contains non-degenerate $k$-APs, then $\aw(G,k) = n + 1$.  For a graph $G$, if $\aw(G,k) = r$, then an \emph{extremal coloring} is an exact ($r-1$)-coloring of $G$ that avoids rainbow $3$-APs.

Notice that at least $k$ colors are needed to have a rainbow $k$-AP. This paper also includes the convention that since a graph cannot be colored with more colors than it has vertices the anti-van der Waerden number of a graph is bounded above by one more than its order.  In the case that $k \ge |G| +1$, then  $\aw(G,k) = |G| +1$.  This is formally stated in Observation \ref{obs1}.

\begin{obs}\label{obs1} 
If $G$ is a graph on $n$ vertices, then $k\leq \aw(G,k) \leq n+1$. If $k \ge n +1$, then $\aw(G,k) = n+1$.
\end{obs}

 In Section \ref{funtools}, results that will be used throughout the paper are established. In Section \ref{p2pnandp3pn}, results are established on $P_m\square P_n$ where $m =2$ or $m = 3$.  In Section \ref{genprod}, these results are used to prove Conjecture \ref{conj1} from a paper authored by Schulte, et. al.
 
 \begin{conj}[\cite{SWY}]\label{conj1}
	If $G$ and $H$ are connected graphs, then
	
		$$\aw(G\square H,3) \le 4.$$
	
\end{conj}

The result from Conjecture \ref{conj1} is used in Section \ref{appsec} to find the anti-van der Waerden number of $P_m\square P_n$ for all $m$ and $n$. 

\section{Fundamental Tools}\label{funtools}

In this section, preliminary results are established that are applicable throughout the remainder of the paper.   A subgraph $H$ of $G$ is \emph{ isometric} if for all $u,v\in V(H)$ $\dist_H(u,v) = \dist_G(u,v)$.

\begin{lem}\label{isopreserveskap}
If $H$ is an isometric subgraph of $G$, then a $k$-AP in $H$ is a $k$-AP in $G$.  If there exists a $k$-AP in $G$ that only contains vertices of $H$, then it is also a $k$-AP in $H$.
\end{lem}

\bpf
Let $\{x_1, x_2, \ldots, x_k\}$ be a $k$-AP in $H$. Since this is a $k$-AP, then $\dist_H(x_i,x_{i+1})=d$ for $1\leq i\leq k-1$. By the definition of isometric subgraph, $\dist_H(x_i,x_{i+1})=\dist_G(x_i,x_{i+1})$.  Hence, $\{x_1, x_2, \ldots, x_k\}$ is a $k$-AP in $G$.\\
Now suppose $\{x_1, x_2, \ldots, x_k\}$ is a $k$-AP in $G$ and $x_i\in V(H)$ for $1\leq i\leq k$. Since $\{x_1, x_2, \ldots, x_k\}$ is a $k$-AP in $G$, then $\dist_G(x_i,x_{i+1})=d'$ for $1\leq i\leq k-1$. Since $H$ is an isometric subgraph, $\dist_G(x_i,x_{i+1})=\dist_H(x_i,x_{i+1})$ for all $1 \leq i\leq k-1$, and therefore, $\{x_1, x_2, \ldots, x_k\}$ is a $k$-AP in $H$.
\epf

\begin{prop}\label{isometric}
	If $H$ is an isometric subgraph of $G$ and $c$ is an exact $r$-coloring of $G$ that avoids rainbow $k$-APs, then $H$ contains at most $\aw(H,k) - 1$ colors.
\end{prop}

\bpf 
Suppose by way of contradiction, $|c(H)|\geq\aw(H,k)$.  This implies $H$ has a rainbow $k$-AP, namely $\{x_1,x_2,\dots,x_k\}$, since every $\aw(H,k)$-coloring of $H$ must have a rainbow $k$-AP by definition. By Lemma \ref{isopreserveskap}, $\{x_1, x_2, \ldots, x_k\}$ is also a $k$-AP in $G$, a contradiction. Hence, any isometric subgraph $H$ of $G$ has at most $\aw(H,k)-1$ colors.
\epf

	Note that Proposition \ref{isometric} ensures that whenever there exists a rainbow $3$-AP in an isometric subgraph of $G$, there is a corresponding rainbow $3$-AP in $G$.  This fact is used frequently without citation in the remainder of this paper.

\begin{lem}\label{blocklemma} 
Let $G = P_m\square P_n$ and $c$ be an exact $r$-coloring of $G$ with $r\ge 3$ that avoids rainbow $3$-APs. If $c(v_{i,j}) = red$ and $c(v_{i-1,j+1}) = blue$, then $c(v_{k,\ell}) \in\{red,blue\}$ when $k\ge i$ and $\ell \ge j+1$ or $k \le i-1$ and $\ell \le j$.  Further, if $c(v_{i,j}) = red$ and $c(v_{i-1,j-1}) = blue$, then $c(v_{k,\ell}) \in \{red,blue\}$ when $k\ge i$ and $\ell \le j-1$ or $k \le i-1$ and $\ell \ge j$.
\end{lem}

\begin{figure}[H]
 \centering
 \begin{tikzpicture}[rotate = 90, scale = .75]
 
    \node (1) at (0,0) [mystyle] {};
    \node (2) at (1,0) [mystyle,fill=blue!40] {B};
    \node (3) at (0,1) [mystyle,fill=red] {R};
    \node (4) at (1,1) [mystyle] {};
    \node (5) at (2,1) [mystyle] {};
    \node (6) at (2,0) [mystyle] {};
    \node (7) at (0,-1) [mystyle] {};
    \node (8) at (1,-1) [mystyle] {};
    \node (9) at (-1,0) [mystyle] {};
    \node (10) at (-1,1) [mystyle] {};
    \node (11) at (0,2) [mystyle] {};
    \node (12) at (1,2) [mystyle] {};
    
    \draw[dashed] (1.4,2) -- (3,2);
    \draw[dashed] (2.4,1) -- (4,1);
    \draw[dashed] (2.4,0) -- (4,0);
    \draw[dashed] (1.4,-1) -- (3,-1);
    \draw[dashed] (-.4,2) -- (-2,2);
    \draw[dashed] (-1.4,1) -- (-3,1);
    \draw[dashed] (-1.4,0) -- (-3,0);
    \draw[dashed] (-.4,-1) -- (-2,-1);
    
    \draw[dashed] (2,-.4) -- (2,-2);
    \draw[dashed] (1,-1.4) -- (1,-3);
    \draw[dashed] (0,-1.4) -- (0,-3);
    \draw[dashed] (-1,-.4) -- (-1,-2);
    \draw[dashed] (2,3) -- (2,1.4);
    \draw[dashed] (1,4) -- (1,2.4);
    \draw[dashed] (0,4) -- (0,2.4);
    \draw[dashed] (-1,3) -- (-1,1.4);
    
    \draw[draw=black] (0.52,0.52) rectangle (4.5,4.5);
    \draw[draw=black] (0.48,0.48) rectangle (-3.5,-3.5);

    \draw (1) -- (2);
    \draw (1) -- (3);
    \draw (3) -- (4);
    \draw (4) -- (2);
    \draw (5) -- (6);
    \draw (2) -- (6);
    \draw (4) -- (5);
    \draw (1) -- (7);
    \draw (7) -- (8);
    \draw (2) -- (8);
    \draw (1) -- (9);
    \draw (9) -- (10);
    \draw (10) -- (3);
    \draw (11) -- (3);
    \draw (11) -- (12);
    \draw (12) -- (4);

\end{tikzpicture}
   \caption{Vertex $R$ (or $v_{i,j}$) is red and vertex $B$ (or $v_{i-1,j+1}$) is blue force the Northwest and Southeast blocks to be red or blue.}
   \label{blockimage}
\end{figure}
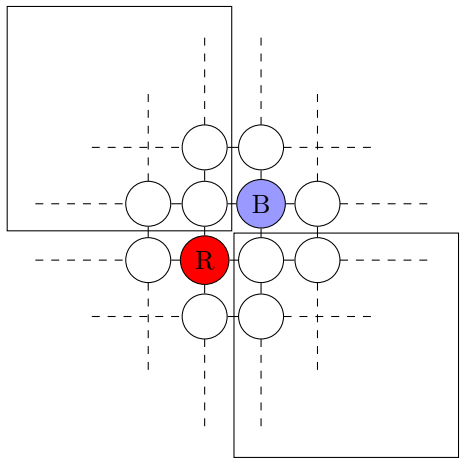

\bpf
Consider the case when $c(v_{i,j}) = red$ and $c(v_{i-1,j+1}) = blue$ (see Figure \ref{blockimage}). Define $v_{k,\ell}$ so that $k\ge i$ and $\ell \ge j+1$. Notice that $\dist(v_{k,\ell},v_{i,j}) = \dist(v_{k,\ell},v_{i-1,j+1}) = k-i+\ell-j$. This means $\{v_{i,j},v_{k,\ell},v_{i-1,j+1}\}$ is a $3$-AP and since $c$ avoids rainbow $3$-APs $c(v_{k,\ell}) \in\{red,blue\}$.  A similar argument can be made in the other three situations.
\epf

\begin{lem}\label{monorow}
	Let $G = P_m\square P_n$ and $c$ be an exact $r$-coloring of $G$ such that $c$ avoids rainbow $3$-APs and  $r\ge 3$.  If $c(v_{i,k}) = \{red\}$ for fixed $i$ and $1 \le k \le n$, $S_1 = \{v_{s,t} \, | \, 1\le s < i,\, 1\le t \le n\}$ and $S_2 =  \{v_{s,t} \, | \, i < s \le m,\, 1 \le t \le n\}$, then $|c(S_i) \cup \{red\}| \le 2$.
\end{lem}

\bpf
	Assume, without loss of generality, that $c(v_{\ell,j}) = blue$ for some $j$ and $i < \ell \le m$ and rows $i+1$ to $\ell -1$ are monochromatic $red$.  By Lemma \ref{blocklemma}, if $c(v_{s,t}) = green$ for $\ell \le s \le m$, $1 \le t \le n$ and $t\neq j$, then either $\{v_{\ell,j}, v_{s,t},v_{\ell-1,j-1}\}$ or $\{v_{\ell,j}, v_{s,t}, v_{\ell -1, j+1}\}$ is rainbow.  This implies that for $t\neq j$, $c(v_{s,t}) \in \{red,blue\}$.  However, using Lemma \ref{blocklemma} with $v_{s,j}$, one of $\{v_{s,j}, v_{\ell,j},v_{s-1,j-1}\}$, $\{v_{s,j}, v_{\ell,j},v_{s-1,j+1}\}$, $\{v_{s,j}, v_{\ell-1,j},v_{s-1,j-1}\}$ or\\ $\{v_{s,j}, v_{\ell-1,j},v_{s-1,j+1}\}$ exists and is rainbow.  Thus, no such $v_{s,t}$ is $green$.  A similar argument applies when $1 \le \ell <i$ and rows $\ell +1$ to $i-1$ are monochromatic $red$.
\epf

Lemma \ref{monorow} says that if there is a monochromatic row in some $P_m\square P_n$, then at most one new color can be introduced below the monochromatic row and at most one new color can be introduced above the monochromatic row.  Note that the argument can be easily applied to monochromatic columns.  Corollary \ref{monocol} states this result.

\begin{cor}\label{monocol}
	Let $G = P_m\square P_n$ and $c$ be an exact $r$-coloring of $G$ such that $c$ avoids rainbow $3$-APs and  $r\ge 3$.  If $c(v_{i,k}) = \{red\}$ for fixed $k$ and $1 \le i \le m$, $S_1 = \{v_{s,t} \, | 1\le s \le m,\, 1 \le t < k\}$ and $S_2 =  \{v_{s,t} \, | \, 1 \le s \le m,\, k < t \le n\}$, then $|c(S_i)\cup\{red\}| \le 2$. 
\end{cor}

Lemma \ref{uniquep2podd} will be useful in combination with Proposition \ref{isometric} in determining the anti-van der Waerden number. 

\begin{lem}\label{uniquep2podd}
	If $G = P_2\square P_{2k+1}$ and $k\ge 1$, then there are precisely two exact $3$-colorings of $G$ that avoid rainbow $3$-APs.  
\end{lem}

\bpf
Let $c$ be an exact $3$-coloring of $G$ that avoids rainbow $3$-APs. Without loss of generality, let $c(v_{1,1}) = red$.  If $c(v_{2,1})= red$, then by Corollary \ref{monocol}, $G$ is colored with at most two colors. Thus, $c(v_{2,1}) = blue$.  If $c(v_{1,2}) = green$, then $\{v_{1,2},v_{1,1}, v_{2,1}\}$ is a rainbow $3$-AP.  Now, consider the following cases.\\

\textit{Case 1:} $c(v_{1,2}) = red$

By Lemma \ref{blocklemma}, $c(v_{2,j}) \in \{red,blue\}$ for $2\le j \le 2k+1$.  If $c(v_{2,2})=blue$, then Lemma \ref{blocklemma} forces both the top and bottom rows to be colored $red$ or $blue$ contradicting that $c$ was an exact $3$-coloring. Thus, $c(v_{2,2})$ must be $red$.  By Corollary \ref{monocol}, columns $3$ through $2k+1$ must be $red$ and $green$, but the bottom row is also $red$ and $blue$; thus, $c(v_{2,j}) = red$ for $3\le j \le 2k+1$. This means $c(v_{1,i}) = green$ for some $3\leq i \leq 2k+1$. If $i \neq 2k+1$, then $\{v_{1,i},v_{2,1},v_{2,i+1}\}$ is a rainbow $3$-AP. Thus, for $i < 2k+1$, $c(v_{1,i}) = red$ and $c(v_{1,2k+1}) = green$.  This is an exact $3$-coloring that avoids rainbow $3$-APs.\\

\textit{Case 2:} $c(v_{1,2})=blue$

If $c(v_{2,2}) = green$ then there exists an obvious rainbow $3$-AP.  If $c(v_{2,2}) = blue$ apply an argument similar to Case $1$ and achieve a symmetric coloring.  Consider if $c(v_{2,2}) = red$.  Let $c(v_{i,j}) = green$ such that $j$ is minimal. By Corollary \ref{monocol}, column $j$ cannot be monochromatic since $red$ and $blue$ appear in column $1$.  If $c(v_{1,j}) = red$ and $c(v_{2,j}) = green$, then $c(v_{1,j-1}) \neq green$ by minimality of $j$, $c(v_{1,j-1}) \neq blue$ by the rainbow $3$-AP $\{v_{1,j-1}, v_{2,2}, v_{2,j}\}$ and $c(v_{1,j-1}) \neq red$ by the rainbow $3$-AP $\{v_{1,j-1}, v_{2,1}, v_{2,j},\}$.  If $c(v_{1,j}) = blue$ and $c(v_{2,j}) = green$, a similar argument can be made.  Finally, the symmetry of column $1$ and $2$ demonstrate that $c(v_{1,j}) \neq green$.\\

Therefore, there are two exact $3$-colorings on $G$ that avoid rainbow $3$-APs. 
\epf

\section{Analysis of $P_2\square P_n$ and $P_3 \square P_n$}\label{p2pnandp3pn}

In this section,  results on $P_2\square P_n$ and $P_3 \square P_n$ are established.  These results are used in conjunction with Proposition \ref{isometric} to obtain other results including Theorem \ref{mainthm}. To begin, first consider the smallest non-trivial $P_m\square P_n$.

\begin{obs}\label{p2p2}
$\aw(P_2\square P_2,3) = 3$
\end{obs}

Almost all of the results in this section require an arbitrary coloring of a graph.  Lemma \ref{uniquep2podd} allows  the elimination of one (or more) colors from half of the graph under the right circumstances.
\begin{prop}\label{p2peven}
For every $k\geq 1$, $\aw(P_2 \square P_{2k},3)=3$.
\end{prop}

\bpf
Consider the graph $P_2\square P_{2k}$ and proceed by induction on $k$. First, if $k=1$, then Observation \ref{p2p2} gives $\aw(P_2 \square P_{2k},3)=3$.

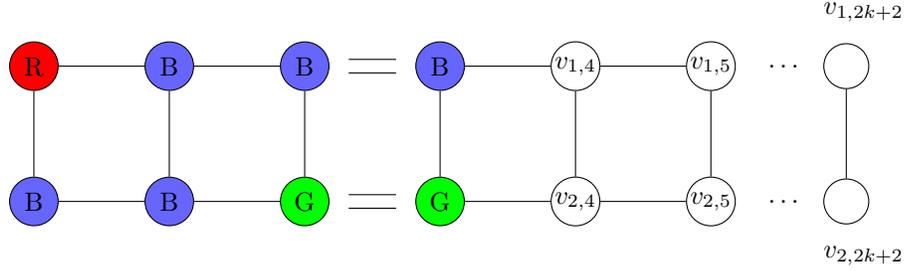
\begin{figure}[h!]
\begin{center}
\begin{tikzpicture}[scale = .9]

  \node (n1) at (1,8)[mystyle, fill=red] {R};
  \node (n2) at (3,8) [mystyle, fill=blue!60] {B};
  \node (n3) at (5,8) [mystyle, fill=blue!60]{B};
  \node (n4) at (7,8) [mystyle, fill=blue!60]{B};
  \node (n5) at (9,8) [mystyle]{$v_{1,4}$};
  \node (n6) at (11,8) [mystyle]{$v_{1,5}$};
  \node (n7) at (13,8) [mystyle]{};
  \node (n9) at (1,6)[mystyle, fill=blue!60]{B};
  \node (n10) at (3,6) [mystyle, fill=blue!60]{B};
  \node (n11) at (5,6) [mystyle, fill=green]{G};
  \node (n12) at (7,6) [mystyle, fill=green]{G};
  \node (n13) at (9,6) [mystyle]{$v_{2,4}$};
  \node (n14) at (11,6) [mystyle] {$v_{2,5}$};
  \node (new) at (12,6) {$\cdots$};
  \node (n15) at (13,6)[mystyle] {};
  \node (n16) at (13,5.2)[mystyle, draw= white] {$v_{2,2k+2}$};
  \node (n17) at (13,8.8)[mystyle, draw= white] {$v_{1,2k+2}$};

  \path (n6) -- node[draw = white, auto=false][fill=white,text width=3mm]{\ldots} (n7);
  \path (n15) -- node[draw = white,auto=false][fill=white,text width=3mm]{\ldots} (n14);

  \foreach \from/\to in {n1/n2,n2/n3,n4/n5,n5/n6,n7/n15,n1/n9,n9/n10,n10/n11,n12/n13,n13/n14,n2/n10,n3/n11,n4/n12,n5/n13,n6/n14}
    \draw (\from) -- (\to);
    \draw (5.65,8.1) -- (6.35,8.1) ;
    \draw (5.65,7.9) -- (6.35,7.9) ;

    \draw (5.65,6.1) -- (6.35,6.1) ;
    \draw (5.65,5.9) -- (6.35,5.9) ;

\end{tikzpicture}\caption{Note double edges indicate vertex identification, so the figure shows $P_2\square P_3 \cup P_2\square P_{2k} = P_2\square P_{2k+2}$.}\label{p2p2kfig}
\end{center}
\end{figure}

For the inductive hypothesis, assume that $\aw(P_2\square P_{2k},3) = 3$. Now consider $P_2\square P_{2k + 2}$ with an exact $3$-coloring $c$ which avoids rainbow $3$-APs. The graph $P_2\square P_{2k+2}$ can be thought of as the union of the two isometric subgraphs formed by columns $1$ through $3$ and columns $3$ through $2k+2$. More technically, let $G_1= P_2\square P_3$ with $V(G_1) = \{v_{1,1},v_{1,2},v_{1,3},v_{2,1},v_{2,2},v_{2,3}\}$ and let $G_2 = P_2\square P_{2k}$ with  $V(G_2) = \{v_{1,3},v_{2,3},v_{1,4},v_{2,4},\ldots , v_{1,2k+2},v_{2,2k+2}\}$. Then $V(G_1) \cap V(G_2) = \{v_{1,3},v_{2,3}\}$ and $G_1 \cup G_2 = P_2\square P_{2k+2}$ (see Figure \ref{p2p2kfig}).  For the following cases, let $c$ be an exact $3$-coloring and let $S = \{v_{1,3},v_{2,3}\}$.\\

Case 1: $|c(S)|=2$.\\
Without loss of generality, let $c(v_{1,3}) = blue$ and $c(v_{2,3}) = red$. By the inductive hypothesis, a third color cannot be introduced into $G_2$ such that there is no rainbow $3$-AP. However, by Lemma \ref{uniquep2podd}, there exists a unique exact $3$-coloring that avoids rainbow $3$-AP's in $G_1$. Without loss of generality, consider the following coloring of $G_1$ where $c(v_{1,1}) = green$, and all other vertices in $G_1$ are colored $blue$ (see Figure \ref{p2p2kfig}).

Now, focusing on the vertex pairs $v_{1,1}$, $v_{2,2}$ and $v_{1,2}$, $v_{2,3}$, Lemma \ref{blocklemma} forces $c(v_{1,j}) = blue$ for $2\le j\le 2k+2$. This however yields the rainbow $3$-AP, $\{v_{1,4},v_{1,1}, v_{2,3}\}$ and this case is complete.\\

Case 2: $|c(S)|=1$\\
Without loss of generality let $c(S) = \{red\}$. By the induction hypothesis and Lemma \ref{uniquep2podd} a maximum of one new color can be added to $G_1$ and one new color can be added to $G_2$ while still avoiding rainbow $3$-APs. Without loss of generality, assume the color introduced in $G_1$ is $blue$ and the color introduced in $G_2$ is $green$. If $c(v_{1,2}) = blue$ then, by Lemma \ref{blocklemma}, $c(v_{1,j})=red$ for $3\le j\le 2k+2$. Now if $c(v_{2,
\ell}) = green$ for some $4\le \ell\le 2k+2$, then Lemma \ref{blocklemma} says that $c(v_{2,1}) = red$, but then $\{v_{2,1},v_{2,\ell},v_{1,2}\}$ is a rainbow $3$-AP.  A similar argument can be made if $c(v_{2,2}) = blue$, so $c(v_{1,2}) = c(v_{2,2}) = red$.\\

Now let $c(v_{1,1})=blue$, then by Lemma \ref{blocklemma} $c(v_{1,j})=red$ for $4\le j\le 2k+2$. This forces $c(v_{2,\ell}) = green$ for some $4\le \ell\le 2k+2$.  If $4\le \ell \le 2k+1$ then $\{v_{2,\ell},v_{1,1},v_{1,\ell+1}\}$ is a rainbow $3$-AP. If $\ell = 2k+2$, then $\{ v_{1,1}, v_{1,k+2}, v_{2,2k+2}\}$ is a rainbow $3$-AP. Similarly, $c(v_{2,1})\neq blue$ which means $|c(G_1)| =1$.  This in turn implies that $|c(G_2)| = 3$ which, as noted earlier, has a rainbow $3$-AP via the inductive hypothesis.
\\

It has been demonstrated that every exact $3$-coloring of $P_2 \square P_{2k+2}$ will result in a rainbow $3$-AP. Thus, $\aw(P_2 \square P_{2k},3)=3$ for all $k\geq 1$.
\epf

\begin{lem}\label{mplusn}
If $G = P_m\square P_n$ and $m+n= 2k+1$ for some $k\geq 1$, then $4 \le \aw(G,3)$.
\end{lem}

\bpf
Consider the exact $3$-coloring $c$ where $c(v_{1,1})=red$, $c(v_{m,n})=blue$ and the remaining vertices are $green$. Note $\dist(v_{1,1}, v_{m,n})=m+n-2$ which, by assumption, is odd so there does not exist a vertex equidistant from both $v_{1,1}$ and $v_{m,n}$, i.e. there is no $3$-AP of the form $\{v_{1,1}, v_{i,j}, v_{m,n}\}$.  This means if a rainbow $3$-AP exists it must be of the form $\{v_{1,1}, v_{m,n}, v_{i,j}\}$ (or similarly $\{v_{m,n}, v_{1,1}, v_{i,j}\}$ which implies there is some vertex $v_{i,j}$ that is distance $m+n-2$ from $v_{1,1}$ or $v_{m,n}$.  However, this cannot happen since $v_{1,1}$ and $v_{m,n}$ are, up to isomorphism, the only two vertices distance $m+n-2$ apart. Thus, an exact $3$-coloring that avoids rainbow $3$-APs has been constructed, therefore $4\leq aw(G,3)$.
\epf

\begin{prop}\label{p2podd}
For every $k\geq 1$, $\aw(P_2 \square P_{2k+1},3)=4$.
\end{prop}

\bpf
 Let $G=P_2 \square P_{2k+1}$. First notice that $4\leq \aw(G,3)$ by Lemma \ref{uniquep2podd}. Now, consider the two isometric subgraphs $G_1 = P_2\square P_2$ and $G_2 = P_2\square P_{2k}$ with $S = V(G_1)\cap V(G_2) = \{ v_{1,2},v_{2,2}\}$. Let $c$ be an exact $4$-coloring of $G$.  Note that $G_1$ and $G_2$ must share at least one color.  If $|c(G_1)| = 2$ and $|c(G_2)| = 2$ then at most three colors have been used.  This implies $|c(G_i)| = 3$ for $i=1$ or $i =2$, but $\aw(G_1,3) = \aw(G_2,3) = 3$ by Observation \ref{p2p2} and Proposition \ref{p2peven}, respectively. Thus, there exists a rainbow $3$-AP in either $G_1$ or $G_2$. Therefore, $\aw(G,3)=4$.
\epf

\begin{prop}\label{p3peven}
For every $k \ge 1$, $aw(P_3 \square P_{2k},3)=4$
\end{prop}

\bpf
Consider the graph $G = P_3 \square P_{2k}$. Since $3+2k = 2(k+1) +1$, then $4\leq \aw(G,3)$ by Lemma \ref{mplusn}.  Let $c$ be an exact $4$-coloring of $G$. Now consider the two isometric subgraphs $G_1$ and $G_2$ each of which are $P_2\square P_{2k}$ graphs where $V(G_1)\cap V(G_2)=\{v_{2,1}, v_{2,2}, \ldots, v_{2,2k}\}$. By Proposition \ref{p2peven}, $G_1$ and $G_2$ must have at most two colors to avoid a rainbow $3$-APs.  This means the coloring $c$ must give a rainbow $3$-AP, thus $\aw(G,3) = 4$.
\epf

\begin{lem}\label{p3p3}
$\aw(P_3\square P_3,3) = 3$
\end{lem}

\bpf
Let $G=P_3\square P_3$ and note that Observation \ref{obs1} gives $3\leq \aw(G,3)$. Let $c$ be an exact $3$-coloring. Consider the two isometric subgraphs $G_1$ and $G_2$ each of which are $P_2\square P_3$ graphs. Let $S = V(G_1)\cap V(G_2) = \{v_{2,1},v_{2,2},v_{2,3}\}$.  If each of these vertices is assigned a different color, then $G$ clearly has a rainbow $3$-AP.\\

Case 1: $|c(S)|=1$.\\
Suppose without loss of generality $c(S)=\{red\}$. By Lemma \ref{monorow}, $G_1$ nor $G_2$ can have three colors.  Without loss of generality, let $c(v_{1,1}) = blue$.  Suppose $c(v_{3,j}) = green$ for $1\leq j\leq 3$. Then, either $\{v_{1,1},v_{2,1},v_{3,1}\}$, $\{v_{3,2},v_{1,1},v_{2,3}\}$ or $\{v_{1,1},v_{2,2},v_{3,3}\}$ is a rainbow $3$-AP.  Therefore, $c(v_{1,1})=red$ and by symmetry 
$$c(v_{1,1})=c(v_{1,3})=c(v_{3,1})=c(v_{3,3})=red.$$  \\
This leaves only $v_{1,2}$ and $v_{3,2}$ uncolored and assigning them the colors $blue$ and $green$ yields the rainbow $3$-AP $\{v_{1,2},v_{2,2},v_{3,2}\}$.\\

Case 2: $|c(S)|=2$.\\
Without loss of generality, let $c(S) = \{blue,green\}$.  A coloring described in Lemma \ref{uniquep2podd} indicates that if a third color is added to $G_1$ or $G_2$, without loss of generality, $c(\{v_{1,2},v_{1,3}, v_{2,1}, v_{2,2}\}) = \{blue\}$, $c(v_{1,1}) = red$,  and $c(v_{2,3})=green$.  If $c(v_{3,1}) = blue$, $c(v_{3,1}) = green$ or $c(v_{3,1}) = red$, then $\{v_{1,1},v_{2,3},v_{3,1}\}$, $\{v_{1,1},v_{2,1},v_{3,1}\}$ or $\{v_{1,2},v_{3,1},v_{2,3}\}$ is a rainbow $3$-AP, respectively.\\

 Therefore, $\aw(G,3) = 3$.\epf

\begin{prop}\label{p3podd}
For every $k \ge 1$, $\aw(P_3\square P_{2k+1},3) = 3$.
\end{prop}

\bpf
Consider the graph $G=P_3\square P_{2k+1}$. First, consider when $k=1$. From Lemma \ref{p3p3}, $\aw(P_3\square P_{2k+1},3) = 3$. Assume that $\aw(P_3\square P_{2k+1},3) = 3$ and now consider the graph $P_3\square P_{2k+3}$.  Recall that $3\leq \aw(P_3\square P_{2k+3},3)$ by Observation \ref{obs1}.  Let $c$ be an exact $3$-coloring of $P_3\square P_{2k+3}$ and consider the two isometric subgraphs $G_1 = P_3\square P_3$ and $G_2 = P_3\square P_{2k+1}$ where $S = V(G_1)\cap V(G_2) = \{v_{1,3},v_{2,3},v_{3,3}\}$. Note that $\aw(G_1,3) = \aw(G_2,3) = 3$ by the base case and induction hypothesis, respectively. Notice that $|c(S)| \neq 2$, otherwise adding a third color to either $G_1$ or $G_2$ would yield a rainbow $3$-AP.  Clearly $|c(S)| \neq 3$, so suppose $|c(S)| = 1$. Without loss of generality, let $c(S) = \{red\}$, $c(V(G_1)) = \{red, blue\}$, and $c(V(G_2)) = \{red, green\}$.\\

If $c(v_{1,2}) = blue$, then $c(v_{1,j}) = red$ for $4\le j \le 2k+3$ by Lemma \ref{blocklemma}. If $c(v_{2,\ell}) = green$ for some $4\leq \ell\leq 2k+2$, then $\{v_{2,\ell},v_{1,2},v_{1,\ell+1}\}$ is a rainbow $3$-AP so $c(v_{2,\ell}) =red$. If $c(v_{2,2k+3})=green$, then $\{v_{1,2},v_{1,k+3},v_{2,2k+3}\}$ is a rainbow $3$-AP. Thus, the color $green$ must only appear in the third row. A similar argument, using $3$-AP $\{v_{3,\ell},v_{1,2},v_{2,\ell+1}\}$, shows that $c(v_{3,\ell}) = red$. If $c(v_{3,2k+3}) = green$, then $c(v_{2,1})$ must be $blue$ since $\{v_{2,1},v_{3,2k+3},v_{1,2}\}$ is a $3$-AP. However, this creates the rainbow $3$-AP $\{v_{3,2k+3},v_{2,1},v_{1,2k+3}\}$. This implies $c(v_{1,2}) = red$ and by symmetry $c(v_{3,2}) = red$.\\

Now, if $c(v_{2,2}) = blue$, then by Lemma \ref{blocklemma}, a third color cannot be introduced in $G_2$. Therefore, all of column two is colored $red$.  \\

If $c(v_{2,1}) = blue$, then by Lemma \ref{blocklemma} a third color cannot be introduced in $G_2$. Thus $c(v_{2,1}) = red$. If $c(v_{1,1}) = blue$, then by Lemma \ref{blocklemma} $c(v_{1,j}) = red$ for $4\leq j\leq 2k+3$. If $c(v_{2,\ell})= green$ for some $4\leq \ell\leq 2k+2$, then $\{v_{1,\ell+1},v_{1,1},v_{2,\ell}\}$ is a rainbow $3$-AP. If $c(v_{2,2k+3}) = green$, then by Lemma \ref{blocklemma} $c(v_{3,1}) = red$ which yields the rainbow $3$-AP $\{v_{1,1},v_{2,2k+3},v_{3,1}\}$. Thus, $c(v_{2,2k+3}) = red$. If $c(v_{3,\ell}) = green$ for $4\leq \ell \leq 2k+2$, then $\{v_{3,\ell},v_{1,1},v_{2,\ell+1}\}$ is a rainbow $3$-AP. Therefore, $c(v_{3,\ell}) = red$. Finally, if $c(v_{3,2k+3}) = green$, then $\{v_{1,1}, v_{2,k+2}, v_{3,2k+3}\}$ is a rainbow $3$-AP. Therefore, any $3$-coloring of $G$ yields a rainbow $3$-AP. 
\epf

\section{General Products}\label{genprod}

In this section, the main result is that the anti-van der Waerden number of Cartesian products of graphs are bounded above by $4$. The section begins with Lemma \ref{diffatmostone} which limits the number of colors that can be introduced in a Cartesian product of graphs.  

\begin{lem}\label{diffatmostone}\cite[Lemma 4.3]{SWY}
Let $G$ be a connected graph on $m$ vertices and $H$ be a connected graph on $n$ vertices.  Let $c$ be an exact $r$-coloring of $G\square H$ with no rainbow $3$-APs. If $G_1, G_2, \dots, G_n$ are the labeled copies of $G$ in $G\square H$, then $|c(V(G_j))\setminus c(V(G_i))|\leq 1$ for all $1\leq i, j \leq n$. 
\end{lem}

It will be useful to find paths in graphs that have at least three colors. This is made possible by the following lemma.

\begin{lem}\label{threevert}
If $G$ is a connected graph on at least three vertices with an exact $r$-coloring $c$ where $r\geq 3$, then there exists a path in $G$ with at least three colors. 
\end{lem}

\bpf 
Choose $u,v \in V(G)$ such that $uv\in E(G)$ and $c(u) = red$ and $c(v) = blue$.  Now, let $w \in V(G)$ such that $\dist(v,w)$ is minimal and $c(w) = green$.  Let $P$ be a shortest path from $v$ to $w$.  If $u$ is on $P$, then $P$ is a path with at least three colors.  If $u$ is not on $P$, let $P'$ be the path from $u$ to $w$ that contains $P$.  Notice $P'$ contains $u$, $v$ and $w$. Therefore, $P'$ has at least three colors.
\epf

\begin{lem}\label{twocolorsbad}
   Assume $G$ and $H$ are connected and consider the graph $G\square H$.  Let $V(H) = \{v_1,v_2,\ldots, v_n\}$, $n \ge 3$, and suppose $c$ is an exact $r$-coloring such that $r\ge 3$, $c$ avoids rainbow $3$-APs and $|c(V(G_i))| \le 2$ for $1\le i\le n$.  If $v_iv_j \in E(H)$, then $|c(V(G_i) \cup V(G_j))| \le 2$.
\end{lem}

\bpf
	If $G_i$ is monochromatic and $G_j$ is monochromatic then the result is immediate.  If $G_i$ is monochromatic and $G_j$ is bichromatic then either $|c(V(G_i) \cup V(G_j))| \le 2$ or $|c(V(G_j))\setminus c(V(G_i))| =2$.  The former is the desired result and the latter contradicts Lemma \ref{diffatmostone}.  Now consider the case where at least one of $G_i$ or $G_j$ has three or more colors.  Without loss of generality, assume $G_i$ has three or more colors.  Then there exists a path with at least three colors, by Lemma \ref{threevert}, in $G_i$.  Let $\rho^{(i)}$ be the shortest path in $G_i$ with three colors and $\rho^{(j)}$ be the corresponding path in $G_j$.  This creates a $P_2\square P_y$ where $y$ is the length of $\rho^{(i)}$.  By Lemma \ref{uniquep2podd} and Proposition \ref{p2peven} there exists a rainbow $3$-AP in $G\square H$.\\
	
	Assume $G_i$ and $G_j$ each have two colors with $|c(V(G_i) \cup V(G_j))| \ge 3$.  Since $|c(V(G_i))\setminus c(V(G_j))| \le 1$, by Lemma \ref{diffatmostone}, then they must share a color.  Without loss of generality, let $c(V(G_i)) = \{red,blue\}$ and $c(V(G_j)) = \{blue,green\}$.  Pick a $red$ vertex, say $v_{i,\alpha}$ , in $G_i$ with a $blue$ neighbor, namely $v$.  Also, choose $v_{j,\beta}$ in $G_j$ such that $c(v_{j,\beta}) = green$.  Let $v_{i,\beta}$ be the vertex in $G_i$ that corresponds to $v_{j,\beta}$ and let $P^{(i)}$ be a shortest path from $v_{i,\alpha}$ to $v_{i,\beta}$ in $G_i$ and $P^{(j)}$ be the corresponding path in $G_j$.  Notice that $P^{(i)}$ and $P^{(j)}$ form an isometric $P_2\square P_x$ in $G\square H$ where $x$ is the length of $P^{(i)}$.  If $P_2\square P_x$ has no $blue$ vertices, then $\{v,v_{i,\alpha},v_{j,\alpha}\}$ is a rainbow $3$-AP.  If $P_2\square P_x$ has a $blue$ vertex and $x$ is even, then there is a rainbow $3$-AP since $\aw(P_2\square P_{2k},3) = 3$ by Proposition \ref{p2peven}.  If $x$ is odd, then by Lemma \ref{uniquep2podd}, so $c(v_{j,\alpha}) =c(v_{i,\beta}) = blue$.  Now, extend to $P_2\square P_x$ to include a corresponding path from $G_k$ where $v_jv_k \in E(H)$, which gives a $P_3\square P_x$ subgraph.  If $P_3\square P_x$ is an isometric subgraph of $G\square H$, then there is a rainbow $3$-AP since $\aw(P_3\square P_{2k+1},3) = 3$, by Proposition \ref{p3podd}.  If $P_3 \square P_x$ is not an isometric subgraph of $G\square H$, then it must correspond to an isometric subgraph $C_3 \square P_x$ of $G\square H$.  Let $v_{k,\beta}$ be the vertex in $G_k$ that corresponds to $v_{j,\beta}$.  However, $c(v_{k,\beta})$ cannot be $red$, $blue$ or $green$ due to $3$-APs $\{v_{i,\beta},v_{j,\beta},v_{k,\beta}\}$, $\{v_{k,\beta},v_{i,\alpha}, v_{j,\beta}\}$ or $\{v_{i,\alpha}, v_{k,\beta}, v_{j,\alpha}\}$. 
\epf

\begin{lem}\label{p2H}
 If $H$ is connected and $|H| \ge 2$, then $\aw(P_2\square H,3) \le 4$.
\end{lem}

\bpf
	Let $c$ be an exact $4$-coloring of $P_2\square H$ with $H_1$ and $H_2$ labeled copies of $H$.  If $|c(V(H_1))| \ge 3$, then, by Lemma \ref{threevert}, there exists a shortest path $P$ with at least three colors in $H_1$.  Then, $P_2 \square P$ is an isometric subgraph of $P_2 \square H$.   By Lemma \ref{uniquep2podd} and Proposition \ref{p2peven} there exists a rainbow $3$-AP in $P_2\square H$  In the case where $|c(V(H_1))| = 1$, then $|c(V(H_2))| \ge 3$, which is  the previous situation.  Finally, consider the case where $|c(V(H_1))| = 2$.  Since $c$ is an exact $4$-coloring of $P_2\square H$, $|c(V(H_1))\backslash c(V(H_2))| = 2$ so by Lemma \ref{diffatmostone} there is a rainbow $3$-AP.
	\epf

The results established thus far come together to show an extremely useful bound on the Cartesian products of graphs in Theorem \ref{mainthm}. This bound demonstrates that the anti-van der Waerden number of any Cartesian product is either $3$ or $4$. 

\begin{thm}\label{mainthm}
If $G$ and $H$ are connected graphs and $|G|,|H| \ge 2$, then $\aw(G\square H,3)\leq 4$.
\end{thm}

\bpf
 If $|H| = |G| = 2$ then $G\square H=P_2 \square P_2$ and by Observation \ref{p2p2} $\aw(G\square H,3)=3\le4$.  Let $c$ be an exact $4$-coloring of $G\square H$ with $V(H) = \{v_1,v_2,\ldots, v_n\}$.   Without loss of generality, assume $|H| \ge 3$.  If $|G| = 2$, then by Lemma \ref{p2H} there is a rainbow $3$-AP.  Now, suppose $|H|,|G| \ge 3$ and define $G_1, G_2, \dots, G_n$ as the labeled copies of $G$ in $G\square H$. Let $\mathcal{P}$ be the path that contains the most colors in some $G_i$, further let it be the shortest such path.  \\

\textit{Case 1:} $\mathcal{P}$ has $3$ or $4$ colors.\\
 Let $\mathcal{P}$ have $x$ vertices and $v_iv_j\in E(H)$.  Also, let $\mathcal{P}'$ be the path in $G_j$ that corresponds to $\mathcal{P}$, note this creates an isometric subgraph $P_2\square P_x$ in $G\square H$.  If $x$ is even, then there is a rainbow $3$-AP since $\aw(P_2\square P_{2k},3) = 3$ for all $k \ge 1$ by Proposition \ref{p2peven} . If $x$ is odd, then  a rainbow $3$-AP is guaranteed by Lemma \ref{uniquep2podd} since path $\mathcal{P}$ has $3$ or $4$ colors.\\

\textit{Case 2:} $\mathcal{P}$ is monochromatic.\\
This implies that each $G_i$ is monochromatic by the definition of $\mathcal{P}$. Since $G\square H$ has $4$ colors, there exists a shortest path $\mathcal{P}'$ in a copy of $H$ that has at least $3$-colors by Lemma \ref{threevert}. However, this is just Case $1$ with the roles of $G$ and $H$ reversed.\\

\textit{Case 3:} $\mathcal{P}$ has two colors.\\
This means that some copy of $G$ has exactly two colors, call this copy $G_d$ and assume the two colors are $red$ and $blue$. By Lemma \ref{diffatmostone}, when the remaining two new colors appear they must both appear either with colors $red$ or $blue$.  Let $yellow$ and $green$ be the two additional colors that are introduced and, without loss of generality, suppose they both appear with $red$.  In particular, let $c(V(G_e)) = \{red,green\}$ and $c(V(G_f))=\{red,yellow\}$.  Now, create an auxiliary coloring $c'$ of  $H$ defined by
\[c'(v_\ell) = \left\{
\begin{array}{cc} 
red & \text{if }c(V(G_\ell)) = \{red\}\\
\mathcal{C} & \text{if }c(V(G_\ell)) = \{\mathcal{C},red\}
\end{array}\right. .
\]

\textit{Subcase 1:} There is no path in $H$, under coloring $c'$, that contains the colors $blue$, $green$ and $yellow$.

Find the smallest subgraph of $H$ that contains $blue$, $green$ and $yellow$, say $c'(v_i) = blue$, $c'(v_j) = green$ and $c'(v_k) = yellow$ and call this smallest subgraph $K$.  This guarantees that $v_i$, $v_j$ and $v_k$ are leaves in the subgraph $K$.  Without loss of generality, assume $\dist(v_i,v_j) \le \dist(v_j,v_k)$. Let $v_{i,\alpha} \in G_i$ such that $c(v_{i,\alpha}) = blue$, $v_{j,\beta} \in G_j$ such that $c(v_{j,\beta}) = green$ and $v_{i,\beta}$ be the vertex in $G_i$ that corresponds to $v_{j,\beta}$.  Let $v_{k,\alpha}$ be the vertex in $G_k$ that corresponds to $v_{i,\alpha}$ and find a shortest path $P$ from $v_{j,\beta}$ to $v_{k,\alpha}$ whose only vertex in $G_j$ is $v_{j,\beta}$.  Now, consider the $3$-AP, $\{v_{i,\alpha}, v_{j,\beta}, u\}$, such that $u$ is a vertex on $P$ since $\dist(v_i,v_j) \le \dist(v_j,v_k)$.  If $c(u) = blue$ or $c(u) = green$ this contradicts the minimality of $K$ or the assumption of the subcase. Therefore, $c(u) \in \{red,yellow\}$ and this $3$-AP is rainbow.\\

\textit{Subcase 2:} There is a path in $H$, under coloring $c'$, that contains $blue$, $green$ and $yellow$.

Let $\mathbb{P}$ be the shortest path in $H$ that contains $blue$, $green$ and $yellow$ and, without loss of generality, assume the path has leaves $v_i$ and $v_k$ with $c'(v_i) = blue$ and $c'(v_k) = yellow$.  Further, assume $v_j$ is the closest $green$ vertex to $v_i$ on $\mathbb{P}$ and $\dist(v_i,v_j) \le \dist(v_j,v_k)$.  Note, there are no other $blue$ or $yellow$ vertices on $\mathbb{P}$, otherwise $\mathbb{P}$ would not be the shortest path that contains $blue$, $green$ and $yellow$.\\

\begin{figure}[H]
  \centering
  \begin{tikzpicture}
    \node (1) at (0,0) [mystyle] {};
    \node (2) at (1,0) [mystyle] {};
    \node (3) at (0,1) [mystyle] {$v_{i,\alpha}$};
    \node (4) at (1,1) [mystyle] {};
    \node (5) at (2,0) [draw=none,fill=none] {$\ldots$};
    \node (6) at (2,1) [draw=none,fill=none] {$\ldots$};
    \node (7) at (3,1) [mystyle] {$v_{i,\beta}$};
    \node (8) at (3,0) [mystyle] {};
    \node (11) at (0,-1) [draw=none,fill=none] {$\vdots$};
    \node (12) at (1,-1) [draw=none,fill=none] {$\vdots$};
    \node (13) at (3,-1) [draw=none,fill=none] {$\vdots$};
    \node (15) at (0,-2) [mystyle] {$v_{j,\alpha}$};
    \node (16) at (1,-2) [mystyle] {};
    \node (17) at (2,-2) [draw=none,fill=none] {$\ldots$};
    \node (19) at (3,-2) [mystyle] {$v_{j,\beta}$};
    \node (20) at (3,-3) [mystyle] {};
    \node (21) at (4,-3) [mystyle] {};
    \node (22) at (4,-2) [mystyle] {};
    \node (23) at (5,-2) [draw=none,fill=none] {$\ldots$};
    \node (23) at (5,-3) [draw=none,fill=none] {$\ldots$};
    \node (24) at (6,-2) [mystyle] {$v_{j,\gamma}$};
    \node (25) at (6,-3) [mystyle] {};
    \node (26) at (3,-4) [draw=none,fill=none] {$\vdots$};
    \node (27) at (4,-4) [draw=none,fill=none] {$\vdots$};
    \node (28) at (6,-4) [draw=none,fill=none] {$\vdots$};
    \node (29) at (3,-5) [mystyle] {$v_{k,\beta}$};
    \node (30) at (4,-5) [mystyle] {};
    \node (31) at (5,-5) [draw=none,fill=none] {$\ldots$};
    \node (32) at (6,-5) [mystyle] {$v_{k,\gamma}$};
    \node (33) at (1.5,2.2) [mystyle, draw = white] {$P$};
    \node (34) at (2,-.5) [mystyle, draw = white] {$P'$};
    \node (35) at (-1.2+8.3,-.5-3) [mystyle, draw = white] {$\rho'$};    
    \node (36) at (-1.2+5.7,-.9) [mystyle, draw = white] {$\rho$};

    \draw (1.5,1) ellipse (2.5cm and .43cm);
    \node (33) at (4.4,1) [draw=none,fill=none] {$G_i$};
    
    \draw (3,-2) ellipse (4cm and .52cm);
    \node (33) at (7.4,-2) [draw=none,fill=none] {$G_j$};
    
    \draw (4.5,-5) ellipse (2.5cm and .43cm);
    \node (33) at (7.4,-5) [draw=none,fill=none] {$G_k$};

    \draw (1) -- (2);
    \draw (1) -- (3);
    \draw (3) -- (4);
    \draw (4) -- (2);
    \draw (7) -- (8);
    \draw (15) -- (16);
    \draw (20) -- (19);
    \draw (20) -- (21);
    \draw (19) -- (22);
    \draw (22) -- (21);
    \draw (24) -- (25);
    \draw (29) -- (30);
    
    \draw [decorate,decoration={brace,amplitude=10pt},xshift=-4pt,yshift=0pt]
(2.7,-2) -- (2.7,1) node [black,midway,xshift=-0.6cm] 
{};

\draw
[decorate,decoration={brace,amplitude=10pt},xshift=-4pt,yshift=0pt]
(0.1,1.5) -- (3.1,1.5) node [black,midway,xshift=-0.6cm] 
{};

    \draw [decorate,decoration={brace,amplitude=10pt},xshift=-4pt,yshift=0pt]
(-.4+7,1-3) -- (-.4+7,-2-3) node [black,midway,xshift=-0.6cm] 
{};

\draw
[decorate,decoration={brace,amplitude=10pt},xshift=-4pt,yshift=0pt]
(0.1+3,1.5-3) -- (3.1+3,1.5-3) node [black,midway,xshift=-0.6cm] 
{};
\end{tikzpicture}
   \caption{Construction of isometric subgraphs of $G\square H$.}\label{subcase2}

 \end{figure}
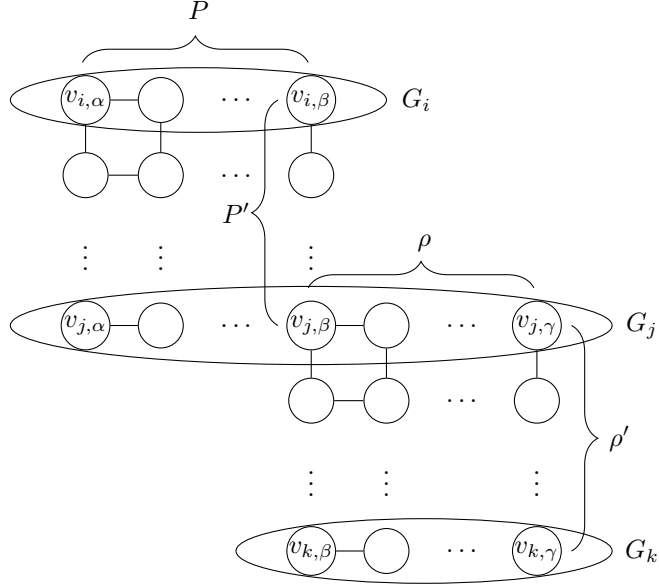

Let $v_{i,\alpha}$ and $v_{j,\beta}$ be in $G_i$ and $G_j$, respectively, so that they are the closest two vertices with $c(v_{i,\alpha})= blue$ and $c(v_{j,\beta})= green$ (see Figure \ref{subcase2} for the following construction).  Let $P$ be a shortest path from $v_{i,\alpha}$ to $v_{i,\beta}$ in $G_i$ and $P'$ be a shortest path from $v_{i,\beta}$ to $v_{j,\beta}$.  Notice that, by minimality of distance from $v_i$ to $v_j$, $P\square P'$ is an isometric subgraph of $G\square H$.  Note that the length of $P'$ is $1$ then there is a rainbow $3$-AP by Lemma \ref{twocolorsbad}.  Assume the length of $P'$ is at least $2$. If $\dist(v_{i,\alpha},v_{j,\beta})$ is even, then there is a $red$ vertex in $P\square P'$, say $u$, such that $\dist(v_{i,\alpha},u) = \dist(u,v_{j,\beta})$ which creates a rainbow $3$-AP.\\

Now, consider the case where $\dist(v_{i,\alpha},v_{j,\beta}) = 2x+1$.  Let $v_{k,\gamma}$ be a vertex in $G_k$ such that $\dist(v_{j,\beta},v_{k,\gamma})$ is minimal and $c(v_{k,\gamma}) = yellow$.  Let $\rho$ be a shortest path from $v_{j,\beta}$ to $v_{j,\gamma}$ in $G_j$ and $\rho'$ be a shortest path from $v_{j,\gamma}$ to $v_{k,\gamma}$, then $\rho \square \rho'$ is an isometric subgraph of $G\square H$.  Note that $c(V(G_{k-1})) = \{red\}$ and $c(V(H_{\gamma - 1})) =\{red\}$ by Lemma \ref{twocolorsbad}.  Define $D_a=\{ v\in V(\rho\square \rho') | \dist(v,v_{j,\beta}) = a\}$ and note that this means $D_0 = \{v_{j,\beta}\}$.  Define $y$ so that $D_y = \{v_{k,\gamma}\}$.  Further, define the distance from $D_s$ to $D_t$ to be $|s-t|$.  If $y < 2x+1$, let $u$ be the vertex on $P'$ or $P$ such that $\dist(u,v_{j,\beta}) = y$.  Then, $c(u) \in \{red,blue\}$ and $\{v_{k,\gamma},v_{j,\beta}, u\}$ is a rainbow $3$-AP.  This means $D_{2x+1} \neq \emptyset$, further, $c(D_{2x+1}) = \{green\}$ because if $v \in D_{2x+1}$, then $\{v_{i,\alpha},v_{j,\beta},v\}$ is a $3$-AP.  This implies that the distance from $D_y$ to either $D_0$ or $D_{2x+1}$ is even.  If $y - 0$ is even, then either 

$$\{v_{k,\gamma}, v_{k-1,\gamma -(y/2 - 1)}, v_{j,\beta}\}\,\, \text{or} \,\,\{v_{k,\gamma}, v_{k-(y/2 - 1)  ,\gamma-1}, v_{j,\beta}\}$$ 

is a rainbow $3$-AP since $c(v_{k-1,\gamma -(y/2 - 1)}) = c(v_{k-(y/2 - 1)  ,\gamma-1}) = red$.  Similarly, if $y - 2x -1$ is even and $z = \frac{y-2x-1}{2}$, then either 

$$\{v_{k,\gamma}, v_{k-1, \gamma - (z-1)}, v_{k -1 - z,\gamma - (z-1)}\}\,\, \text{or}\,\, \{v_{k,\gamma}, v_{k-(z-1),\gamma-1},v_{k-(z-1),\gamma - 1 -z}\}$$

\noindent is a rainbow $3$-AP.\\

Therefore, each case yields a rainbow $3$-AP so $\aw(G\square H,3) \le 4$.
\epf

\section{Application to $P_m\square P_n$}\label{appsec}

In this Sections \ref{funtools} and \ref{p2pnandp3pn} results for $m=2$ and $m=3$ were established.  The result of Theorem \ref{mainthm} is used, with earlier results, to determine $\aw(P_m \square P_n,3)$ for all $m$ and $n$.  It is interesting to note that the pattern for the small values of $m$ does not continue when considering large values of $m$.  Essentially, there are `more' $3$-APs which forces the anti-van der Waerden number to always be $4$. First notice that Lemma \ref{mplusn} and Theorem \ref{mainthm} give the result of Corollary \ref{cor} immediately.

\begin{cor}\label{cor}
If $G = P_m\square P_n$ and $m+n= 2k+1$ for some $k\geq 1$, then $\aw(G,3) = 4$.
\end{cor}

Lemma \ref{mplusneven} gives the final lower bound to determine the anti-van der Waerden number for all $P_m\square P_n$.\\

\begin{lem}\label{mplusneven}
If $m\geq 4$, $n\geq 4$ and $m+n= 2k$ for some $k\geq 1$, then $4\le \aw(P_m\square P_n,3)$. 
\end{lem}

\bpf
Let $G = P_m \square P_n$. Define

$$c(v_{i,j}) = \left\{\begin{array}{cc} red & \text{ if $i=1$ and $j=2$ or $i = 2$ and $j=1$}\\
							blue & \text{if $i = m$ and $j=n$}\\
							green & \text{otherwise}\end{array}\right.$$
	
Note that if a rainbow $3$-AP exists it must contain vertex $v_{m,n}$ and either $v_{1,2}$ or $v_{2,1}$. Let $S=\{v_{m,n},v_{1,2},v_{2,1}\}$.  Note that  $\dist(v_{1,2},v_{m,n}) = \dist(v_{2,1}, v_{m,n})=m+n-3$ which, by assumption, is odd. Therefore, there does not exist a vertex equidistant from $v_{2,1}$ and $v_{m,n}$ or equidistant from $v_{1,2}$ and $v_{m,n}$.  This means a rainbow $3$-AP cannot exist in the order of $\{v_{2,1}, v_{i,j}, v_{m,n}\}$ or $\{v_{1,2},v_{i,j},v_{m,n}\}$.\\

This means any rainbow $3$-AP must exist in the order of $\{v_{m,n}, v_{2,1}, v_{i,j}\}$ or $\{v_{m,n}, v_{1,2}, v_{i,j}\}$ (or the reverse order) where $v_{i,j}\notin S$.  Note that $v_{i,j}$ must be distance $m+n -3$ from one of the vertices in $S$, but the only vertices distance $m+n-3$ from any vertex in $S$ are already in $S$ thus $v_{i,j}$ does not exist.   Therefore, $c$ avoids rainbow $3$-APs so $4\le \aw(G,3)$.
\epf

Using Theorem \ref{mainthm}, Corollary \ref{cor} and Lemma \ref{mplusneven} gives Corollary \ref{mnbig}.

\begin{cor}\label{mnbig}

If $m\geq 4$ and $n\geq 4$, then $\aw(P_m\square P_n,3) = 4$. 

\end{cor}

Finally, combining Propositions \ref{p2peven}, \ref{p2podd}, \ref{p3peven}, \ref{p3podd} and Corollary \ref{mnbig} gives a function to determine $\aw(P_m\square P_n,3)$ for all $m$ and $n$.

\begin{thm}
	For $m,n\ge 2$, 
	$$\aw(P_m\square P_n,3) = \left\{\begin{array}{cc}
		3 & m=2 \,\, \text{and}\,\, $n$\,\, \text{is even or } m= 3\,\, \text{and}\,\, n\,\, \text{is odd}\\
		4 & \text{otherwise}

		\end{array}\right. .$$
\end{thm}

This research was supported by the University of Wisconsin-La Crosse Dean's Distinguished Fellowship.

\end{document}